\documentstyle[11pt,amscd,amssymb]{amsart}

\newcommand{\seq}[3]{{#1}_{#2}, \ldots, {#1}_{#3}}
\newcommand{\Dws}{D^{(w,\Sigma)}}
\newcommand{\Y}{\Sigma \times {{\Bbb S}}^1}
\newcommand{\Spz}{\text{Sp}\, (2g,{\Bbb Z})}
\newcommand{\Ct}{{\Bbb C}[[t]]}
\newcommand{\Ctabg}{{\Bbb C}[[t]][\alpha,\beta,\gamma]}
\newcommand{\Cabg}{{\Bbb C}[\alpha,\beta,\gamma]}
\newcommand{\Ctab}{{\Bbb C}[[t]][\alpha,\beta]}
\newcommand{\Ctrb}{{\Bbb C}[[t]][\bar{\beta}]}
\newcommand{\Crb}{{\Bbb C}[\bar{\beta}]}
\newcommand{\Ctrbg}{{\Bbb C}[[t]][\bar{\beta},\gamma]}
\newcommand{\Crbg}{{\Bbb C}[\bar{\beta},\gamma]}
\newcommand{\ima}{{\bf i}}
\newcommand{\whff}{\widetilde{HFF}{}^*_g}

\newcommand{\rb}{\bar{\beta}}
\newcommand{\rg}{\bar{g}}
\newcommand{\rT}{\bar{T}}
\newcommand{\rJ}{\bar{J}}
\newcommand{\rR}{\bar{R}}

\newcommand{\la}{\langle}
\newcommand{\ra}{\rangle}
\newcommand{\surj}{\twoheadrightarrow}
\newcommand{\inc}{\hookrightarrow}
\newcommand{\ar}{\rightarrow}
\newcommand{\bd}{\partial}
\newcommand{\x}{\times}
\newcommand{\ox}{\otimes}
\newcommand{\iso}{\cong}

\newcommand{\point}{\text{pt}}
\newcommand{\CP}{{\Bbb C \Bbb P}}

\newcommand{\rk}{\text{rk}}

\newcommand{\Gr}{\text{Gr}}

\newcommand{\Sym}{\text{Sym}}

\newcommand{\PD}{\text{P.D.}}

\newcommand{\cF}{{\cal F}}

\newcommand{\cH}{{\cal H}}

\newcommand{\cI}{{\cal I}}

\renewcommand{\AA}{{\Bbb A}}
\newcommand{\CC}{{\Bbb C}}

\newcommand{\SS}{{\Bbb S}}
\newcommand{\ZZ}{{\Bbb Z}}

\renewcommand{\a}{\alpha}
\renewcommand{\b}{\beta}

\newcommand{\g}{\gamma}
\newcommand{\e}{\varepsilon}

\newcommand{\h}{\theta}
\renewcommand{\l}{\lambda}

\newcommand{\p}{\phi}
\newcommand{\q}{\psi}

\renewcommand{\S}{\Sigma}
\newcommand{\D}{\Delta}
\renewcommand{\L}{\Lambda}

\newcommand{\frs}{{\frak s}}

\newcommand{\frm}{{\frak m}}

\theoremstyle{plain}
\begingroup
\theorembodyfont{\sl}
\newtheorem{thm}{Theorem}[section]
\newtheorem{cor}[thm]{Corollary}
\newtheorem{lem}[thm]{Lemma}
\newtheorem{prop}[thm]{Proposition}

\endgroup

\theoremstyle{definition}
\newtheorem{defn}[thm]{Definition}

\theoremstyle{remark}
\newtheorem{rem}[thm]{Remark}

\title{Higher type adjunction inequalities for Donaldson invariants}

\author{Vicente Mu\~noz}
\address{Departamento de \'Algebra Geometr\'{\i}a y Topolog\'{\i}a \\ 
  Facultad de Ciencias \\ Universidad de M\'alaga \\ 29071 M\'alaga
\\ Spain}
\email{vmunoz@@agt.cie.uma.es}
\thanks{\hbox{$^*$}Partially supported by DGES through Spanish
Research 
  Project PB97-1095. 
  }
\date{December, 1998. Revised May, 1999.}

\begin{document}

\maketitle

\begin{abstract}
  We prove new adjunction inequalities for embedded surfaces 
  in four-manifolds with non-negative 
  self-intersection number using the Donaldson
  invariants. These formulas are completely analogous to the ones
  obtained by
  Ozsv\'ath and Szab\'o~\cite{OS} using the Seiberg-Witten
  invariants. To 
  prove these relations, we give a fairly explicit description of the 
  structure of the Fukaya-Floer homology of a surface times a circle.
  As an aside, we also relate the Floer homology of a surface times a 
  circle with the cohomology of some symmetric products of the
surface.
\end{abstract}

\section{Introduction}
\label{sec:intro}

In this paper, we prove new adjunction inequalities for embedded
surfaces of 
non-negative self-intersection in a four-manifold $X$ by making use
of the 
Donaldson invariants of $X$. Our results are parallel to those of
Ozsv\'ath and Szab\'o~\cite{OS} on adjunction inequalities from 
Seiberg-Witten theory, and once more, they provide a confirmation of
the (yet conjectural) equivalence between the Donaldson invariants
and the 
Seiberg-Witten invariants~\cite{W} and give a guideline on how this
correspondence 
should be established. The scope of application of these new
adjunction
inequalities is somewhat narrow at present, as they provide new
information
only for $4$-manifolds not of simple type or with $b_1 >0$. So far no
example
of $4$-manifold not of simple type with $b^+>1$ has been given. This
same 
situation applies to the results in~\cite{OS}. In turn the
description and use
of the Fukaya-Floer homology is very enlightening and clarifies the
meaning of
the basic classes~\cite{basic} and of the finite type condition.

The input for the analysis carried out here is the Floer
homology~\cite{hf}
and Fukaya-Floer homology~\cite{hff} of the three-manifold $Y=\Y$, 
where $\S$ is a surface of genus $g\geq 1$. 
Using the descriptions, in terms of generators and relations, of (the
ring
structure of) $HF^*(Y)$ and $HFF^*(Y,\SS^1)$ given in~\cite{hf}
and~\cite{hff},
we are able to analyse the pieces of their artinian decompositions
and check that 
they look like (a deformation of) the cohomology of some symmetric
products of 
the surface $\S$. This turns out to provide enough information to get
new
relations in the effective
(Fukaya-)Floer homology, which are recast in the shape
of adjunction inequalities for $4$-manifolds with $b^+>1$.

In order to state our main results, let us set up some notation.
Donaldson invariants for a (smooth, compact, oriented) $4$-manifold
$X$
with $b^+>1$ (and with a homology orientation) are defined as linear
functionals~\cite{KM}
$$
  D^w_X: \AA(X)= \Sym^*(H_0(X) \oplus H_2(X)) \ox \L^* H_1(X) \ar
\CC,
$$
where $w \in H^2(X;\ZZ)$. As for the grading of
$\AA(X)$, we give degree $4-i$ to the elements in $H_i(X)$
(complex coefficients are understood for the homology $H_*(X)$ and 
cohomology $H^*(X)$).
We shall denote by $x \in H_0(X;\ZZ)$ the class of a point.
Recall~\cite{hff} that we say that a $4$-manifold $X$ with $b^+>1$ is
of
finite type if there exists $n \geq 0$ such that $D_X^w((x^2-4)^n
z)=0$,
for all $z \in \AA(X)$, and that we call the minimum such $n$ the
order
of finite type of $X$. All $4$-manifolds $X$ with $b^+>1$ are of
finite 
type~\cite{hff} and the order of finite type does not depend on 
$w\in H^2(X;\ZZ)$ (see~\cite{basic}). 

Our first main result is a bound on the order of finite type
for $4$-manifolds, which greatly improves the one given
in~\cite[proposition
7.3]{hff} in the case where $X$ has $b_1 >0$.

\begin{thm}
\label{thm:0}
  Let $X$ be a $4$-manifold with $b^+>1$
  and suppose that there is an embedded surface
  $\S \subset X$ of genus $g$ and either with self-intersection
$\S^2>0$ 
  or with $\S^2 =0$ and $\S$ representing an odd homology class. Then
the order
  of finite type of $X$ is less than or equal to $g$.
\end{thm}

In~\cite{basic} we saw that the Donaldson invariants of a
$4$-manifold
$X$ with $b^+>1$ always produce a distinguished set of cohomology
classes
$K_i$ (called basic classes), regardless of $X$ being of simple type
or
having $b_1=0$, extending thus the notion of basic classes given 
in~\cite{KM}. These were obtained by looking at the asymptotics of
$D^w_X$
for large degrees. Also the set of basic classes is independent of $w
\in H^2(X;\ZZ)$.
Recall the following criterium extracted from~\cite[proposition
15]{basic}.

\begin{prop}
\label{prop:basic}
  Let $X$ be a $4$-manifold with $b^+>1$ and $w \in H^2(X;\ZZ)$. 
  Then $K \in H^2(X;\ZZ)$ is a basic class for $X$ if and only if
there
  exists $z\in \AA(X)$ such that
  $D^w_X(z e^{tD+\l x})=e^{Q(tD)/2+2\l+K \cdot tD}$, for all $D \in
H_2(X)$.
  \hfill $\Box$
\end{prop}

We want to define the order of finite type for a single basic class
$K$
for $X$. This should be the analogue of the expected
dimension $d(\frs)$ of the Seiberg-Witten moduli space associated to
a 
spin$^{\CC}$ structure $\frs$ with $c_1(\frs)=K$ and non-zero
Seiberg-Witten invariant
$SW_{X,\frs}$ (see~\cite{OS}). In that case there exists $SW_{X,\frs}
(x^a\g_{1}\cdots \g_{r}) \neq 0$ with $2a+r=d(\frs)$, and $\g_i \in
H_1(X)$. The
heuristic comparison of the Donaldson and Seiberg-Witten theories
makes correspond
$x$ in Seiberg-Witten theory to $\wp=(x^2-4) \in \AA(X)$ in Donaldson
theory. So we define
$$
  \tilde{\AA}(X)= \CC[ \wp ] \ox \L^* H_1(X),
$$
where we give the following grading: $d(\wp)=2$, and for any
$\g \in H_1(X)$, $d(\g)=1$. Note that this grading has nothing to do 
with the grading of $\AA(X)$ 
under the inclusion $\tilde{\AA}(X) \subset \AA(X)$. The fact that
any $4$-manifold $X$ with $b^+>1$ is of finite type is transcribed as
$D^w_X(\wp^n z)=0$, for any $z \in \AA(X)$ and sufficiently large
$n$.
Now we can give the following two definitions.

\begin{defn}
\label{def:b}
  Let $X$ be a $4$-manifold with $b^+>1$, and let $b \in
\tilde{\AA}(X)$.
  If $K \in H^2(X;\ZZ)$ then we say that $K$ is a basic class for 
  $D^w_X(b\,\bullet)$ if there exists $z\in \AA(X)$ such that
  $D^w_X(b z e^{tD+\l x})=e^{Q(tD)/2+2\l+K \cdot tD}$, for all 
  $D \in H_2(X)$. This condition is independent of $w\in H^2(X;\ZZ)$.
\end{defn}

Using~\cite[theorem 6]{basic} we have that for $b \in \tilde{\AA}(X)$
and any 
homogeneous $z \in \L^* H_1(X)$, 
there exist polynomials $p_{i,z}, q_{i,z} \in \Sym^* H^2(X)\ox
\CC[\l]$, for every
basic class $K_i$, such that
\begin{equation}
\label{eqn:h}
  D^w_X(b z \, e^{tD+\l x})=\sum p_{i,z}(tD,\l)
e^{Q(tD)/2+2\l+K_i\cdot tD} 
  +\sum q_{i,z}(tD,\l) e^{-Q(tD)/2-2\l+\ima K_i\cdot tD},
\end{equation}
for all $D \in H_2(X)$. Then the condition that $K$ be a basic class
for
$D^w_X(b\,\bullet)$ is equivalent to the existence of some 
homogeneous $z\in \L^* H_1(X)$
such that $p_{i,z} \neq 0$ in~\eqref{eqn:h} for $K=K_i$.

\begin{defn}
\label{def:d(K)}
  Let $X$ be a $4$-manifold with $b^+>1$, and  $K \in H^2(X;\ZZ)$ a
basic
  class for $X$. Then we define the order of finite type of $K$ to be
  $$ 
 d(K)=\text{max}\{d(b)|\text{$K$ is a basic class for
$D^w_X(b\,\bullet)$}, \,
  b \in \tilde{\AA}(X) \}.
  $$
\end{defn}

\begin{rem}
\label{rem:rm}
  We leave the proof of the following characterization of $d(K)$ to
the reader.
  Collecting together (as $z$ runs through an homogeneous basis of
$\L^* H_1(X)$)
  all the polynomials from~\cite[theorem 6]{basic}, we have
  $P_i, Q_i \in \Sym^* H^2(X)\ox \CC[\l] \ox \L^* H^1(X)$, for every
  basic class $K_i$, such that
$$
  D^w_X(z \, e^{tD+\l x})=\sum P_i(tD,\l,z) e^{Q(tD)/2+2\l+K_i\cdot
tD} 
  +\sum Q_i(tD,\l,z) e^{-Q(tD)/2-2\l+\ima K_i\cdot tD},
$$
  for all $D \in H_2(X)$ and $z \in \L^* H_1(X)$. Then $d(K_i)= \deg
P_i=\deg Q_i$, where
  the elements in $H^2(X)$ have degree $2$, $\l$ has degree $2$ and
the elements in
  $H^1(X)$ have degree $1$.
\end{rem}

The main results that we prove are refinements, for four-manifolds
not of simple type or with $b_1 >0$, of the adjunction 
inequality proved by Kronheimer and Mrowka~\cite{KM}.
They are analogues to theorem 1.1,
theorem 1.3 and theorem 1.4 of~\cite{OS}, which are established in
the context of 
Seiberg-Witten invariants.

\begin{thm}
\label{thm:A}
  Let $X$ be a $4$-manifold with $b^+>1$ and let $\S \subset X$ be an
embedded
  surface of genus $g\geq 1$ either with self-intersection $\S^2>0$
or with 
  $\S^2 =0$ and $\S$ representing an odd homology class. 
  Let $b \in \tilde{\AA}(\S)$.
  If $K$ is a basic class for $D^w_X(b\,\bullet)$ then 
  $$
   |K \cdot \S| + \S^2 + d(b) \leq 2g-2.
  $$
\end{thm}

Note in particular that theorem~\ref{thm:0} follows from
theorem~\ref{thm:A}.
With addition of the condition $b_1=0$ on $X$ we get the stronger
result
\begin{thm}
\label{thm:B}
  Let $X$ be a $4$-manifold with $b^+>1$ and $b_1=0$, 
  and let $\S \subset X$ be an embedded surface 
  of genus $g\geq 1$ either with self-intersection $\S^2>0$ or with 
  $\S^2 =0$ and $\S$ representing an odd homology class. 
  If $K$ is a basic class for $X$ then we have the following
adjunction
  inequality 
  $$
   |K \cdot \S| + \S^2 + 2 d(K) \leq 2g-2.
  $$
\end{thm}

As in~\cite{OS}, theorem~\ref{thm:B} may be generalised as
\begin{thm}
\label{thm:C}
  Let $X$ be a $4$-manifold with $b^+>1$ 
  and let $\imath:\S \inc X$ be an embedded surface 
  of genus $g\geq 1$ either with self-intersection $\S^2>0$ or with 
  $\S^2 =0$ and $\S$ representing an odd homology class. Let $l$ be
an integer
  so that there is a symplectic basis $\{\g_i\}_{i=1}^{2g}$ of
$H_1(\S)$ with 
  $\g_i\cdot\g_{g+i}=1$, $1\leq i \leq g$, satisfying that 
  $\imath_* (\g_j)=0$ in $H_1(X)$ for $j=1,\ldots, l$.  
  Let $b \in \tilde{\AA}(\S)$ be an element of degree $d(b) \leq
l+1$.
  If $K$ is a basic class for $D_X^w(b\, \bullet)$ then we have 
  $$
   |K \cdot \S| + \S^2 + 2 d(b) \leq 2g-2.
  $$
\end{thm}

A simple application is the following. 
Take $X$ and $\S$ as in the statement of 
theorem~\ref{thm:A} such that there is a basic class $K$ of $X$ with
$|K \cdot \S| + \S^2 = 2g-2$ (for instance the manifolds $B_g$ or
$C_g$ 
of~\cite[definition 25]{genusg} will do, for any $g \geq 1$). Let $l
\geq 1$ 
and consider the connected sum $X'=X \# l (\SS^1 \x \SS^3)$. 
Put $\g_i$ for the loop corresponding to the $\SS^1$ factor in the
$i$-th
copy of $\SS^1\x\SS^3$ and let $T_i \subset \SS^1 \x\SS^3$ be the 
(homologically trivial) torus of the natural elliptic fibration,
for $1 \leq i \leq l$.
Let $\S' =\S \# T_1 \# \cdots \# T_l \subset X'$ be the genus $g+l$ 
surface obtained by performing internal
connected sums. As $D^w_{X'} (\g_1\cdots \g_l z)= D^w_X (z)$, for any
$z \in \AA(X)$,
theorem~\ref{thm:C} is sharp for $X'$ and $\S'$ with
$b=\g_1\cdots \g_l$. So the homology class of $\S'$ cannot be
represented by 
an 
embedded surface $S \subset X'$ with $H_1(T_1) \oplus \cdots \oplus
H_1(T_l) \subset 
H_1(S)$ and genus strictly less than $g+l$.

In section~\ref{sec:hff} we describe the Fukaya-Floer homology
$HFF^*_g$ 
of the three manifold $Y=\Y$, where $\S$ is a surface of genus $g\geq
1$.
This is a finite graded commutative $\Ct$-algebra generated by
elements 
$\a$, $\b$ and $\q_i$, $1\leq i\leq 2g$, canonically associated to
generators 
of the homology $H_*(\S)$. The eigenvalues of the $\Ct$-linear 
(commuting) endomorphisms 
given by multiplication by $\a$, $\b$ and $\q_i$, $1\leq i\leq 2g$,
form a 
discrete set. Actually there is a Jordan decomposition of $HFF^*_g$
with
respect to all of these endomorphisms, which corresponds with the 
artinian decomposition of $HFF^*_g$ interpreted as a module over the
free
algebra $\Ct[\a,\b]\ox \L^* (\seq{\q}{1}{2g})$.
Each of the pieces $\cH_r$ of the decomposition $HFF^*_g=\oplus
\cH_r$
is labelled by an integer $r$ running
between $-(g-1)$ and $g-1$, and is the localization of $HFF^*_g$ at 
a corresponding prime ideal. The piece $\cH_r$ controls the basic
classes $K$ such that $K\cdot \S=2r$, 
for any $4$-manifold $X$ with $b^+>1$ and an embedded
surface $\S\subset X$ with $\S^2=0$ (and representing an odd element
in homology).

In section~\ref{sec:relations} 
we get the new relations coming out of the description of the
Fukaya-Floer 
homology of $Y$ and also how imposing ``extra'' conditions produces
more relations.
This is used in section~\ref{sec:proofs} in a fairly straightforward
way
to prove our main results: theorem~\ref{thm:0} and
theorems~\ref{thm:A}--\ref{thm:C}.

For completeness, in section~\ref{sec:hf} we show how the artinian 
decomposition of the Floer homology $HF^*(Y)$ is related to the
cohomology
of some symmetric products of $\S$, looking like a deformation of
their natural 
ring structures. This is what one would expect it to be, as 
the Seiberg-Witten-Floer homologies of $Y$ (labelled by the
spin$^{\Bbb C}$ structures 
on $Y$) are a deformation of the cohomology rings
of symmetric products of $\S$, and probably isomorphic to the quantum
cohomology of such
spaces (see~\cite{wang-yo}).

\section{Structure of the Fukaya-Floer homology of $\Y$}
\label{sec:hff}

Fix a surface $\S$ of genus $g \geq 1$ and consider
the $3$-manifold $Y=\Y$ with the $SO(3)$-bundle
with $w_2=\PD[\SS^1] \in H^2(Y;\ZZ/2\ZZ)$ and the loop 
$\SS^1 \subset Y=\Y$. Associated to this triple we have defined
Fukaya-Floer (co)homology groups $HFF^*_g=
HFF^*(Y,\SS^1)$ which were (partially) determined in~\cite[section
5]{hff}.
In this section it is our purpose to describe
the artinian decomposition of the ring structure of $HFF^*_g$.

There are two basic properties of the Fukaya-Floer homology. The
first
one is the existence of relative Donaldson invariants.
For every $4$-manifold $X_1$ with
boundary $\bd X_1=Y$, $w_1 \in H^2(X_1;\ZZ)$ such that $w_1|_Y=w \in
H^2(Y;\ZZ/2\ZZ)$, $z \in \AA(X_1)$ and $D_1 \subset X_1$ a $2$-cycle
with $\bd D_1=\SS^1$, one has a relative invariant
$$
  \p^w(X_1,z e^{tD_1}) \in HFF^*_g.
$$
The second basic property of the Fukaya-Floer homology is the 
existence of a pairing
$$
   \la\/,\ra:HFF^*_g \ox HFF^*_g \ar \Ct,
$$
satisfying a gluing property
for the Donaldson invariants~\cite[theorem 3.1]{hff}.
Let $X=X_1 \cup_Y X_2$ be a $4$-manifold, split in two
$4$-manifolds $X_1$ and $X_2$ with boundary $\bd X_1=-\bd X_2 =Y$,
and $w \in H^2(X;\ZZ)$ satisfying 
$w \cdot \S \equiv 1\pmod 2$. Put $w_i =w|_{X_i} \in H^2(X_i;\ZZ)$. 
Let $D \in H_2(X)$ be decomposed as $D=D_1 +D_2$ with $D_i \subset
X_i$, 
$i=1,2$, $2$-cycles with $\bd D_1=-\bd D_2=\SS^1$. For $z_i \in 
\AA(X_i)$, $i=1,2$, it is
  $$
     \Dws_X(z_1z_2e^{tD})=
     \la\p^{w_1}(X_1,z_1e^{tD_1}),\p^{w_2}(X_2,z_2e^{tD_2})\ra,
  $$
where $\Dws_X=D^w_X+D^{w+\S}_X$. If $X$ has $b^+=1$ then
the invariants are calculated for metrics on $X$ giving a long neck
to the splitting $X=X_1 \cup_Y X_2$.

Let $A=\S \x D^2$ be the product of $\S$ times a $2$-dimensional disc
and
consider the horizontal section $\D=\point\x D^2 \subset A$. Put
$w=\PD [\D]
\in H^2(A;\ZZ)$. Let $\{\g_i\}$ be a symplectic basis for $H_1(\S)$
with $\g_i \cdot \g_{g+i}=1$, $1\leq i \leq g$.
We have the following elements
\begin{equation}
   \left\{ \begin{array}{l} \a= 2 \, \p^w(A,\S\,e^{t\D}) \in HFF^2_g
   \\ \q_i= \p^w(A,\g_i\,e^{t\D}) \in HFF^3_g, \qquad 1\leq i \leq 2g
   \\ \b= - 4 \, \p^w(A,x\,e^{t\D}) \in HFF^4_g
    \end{array} \right.
\label{eqn:gen}
\end{equation}
The Fukaya-Floer homology $HFF^*_g$
is a $\Ct$-algebra~\cite[section 5]{hff} generated by
the elements~\eqref{eqn:gen} and the product for $HFF_g^*$ is
actually
determined by the property 
$$
  \p^w(A, z_1e^{t\D})\p^w(A, z_2 e^{t\D})= \p^w(A, z_1z_2 e^{t\D}),
$$
for $z_1, z_2\in \AA(\S)$. 
The mapping class group of $\S$ acts on $HFF^*_g$ factoring through
an
action of the symplectic group $\Spz$ on $\{\q_i\}$.
The invariant part, $(HFF^*_g)_I$, is generated by $\a$, $\b$
and $\g=-2 \sum \q_i\q_{g+i}$.
For $0 \leq k \leq g$, we define the primitive component of 
$\L^k =\L^k (\seq{\q}{1}{2g})$ as
$$
  \L_0^k=\L_0^k (\seq{\q}{1}{2g})=\ker (\g^{g-k+1}:\L^k \ar
\L^{2g-k+2}).
$$
The spaces $\L^k_0$ are irreducible $\Spz$-representations. We
have the following structural result

\begin{thm} {\rm (\cite[theorem 5.3]{hff})}
\label{thm:hff.struct}
  Let $\S$ be a surface of genus $g\geq 1$. Then
  $HFF_g^*$ is, as $\Spz$ representation,
\begin{equation}
   HFF^*_g= \bigoplus_{k=0}^{g-1} \L_0^k  \ox {\Ctabg \over J_{g-k}},
\label{eqn:hff.struct}
\end{equation}
  where $J_r=(R^1_r, R^2_r,R^3_r)$ and $R^i_r$ are defined 
  recursively by setting 
  $R^1_0=1$, $R^2_0=0$, $R^3_0=0$ and putting, for all $0\leq r \leq
g-1$,
  $$
   \left\{ \begin{array}{l} R_{r+1}^1 = (\a+f_{11}) R_r^1 + r^2 
     (1+f_{12}) R_r^2 +f_{13} R_r^3
   \\ R_{r+1}^2=(\b+(-1)^{r+1}8+f_{21})R_r^1+f_{22}R_r^2+ ({2r \over
r+1}+
     f_{23})R_r^3 
   \\ R_{r+1}^3 = \g  R_r^1
    \end{array} \right.
  $$
  for some (unknown) functions $f_{ij} \in t\Ct[\a,\b,\g]$,
  dependent on $r$ and $g$. \hfill $\Box$
\end{thm}

The actual shape of the relations is not of much importance for our
purposes.
We intend now to understand the artinian decomposition of $HFF^*_g$
in
a fairly clean way to turn the algebra information of $HFF^*_g$ into
properties for the Donaldson invariants. To shorten the notation, we
shall
write
 \begin{equation}
   T_{g,k}= {\Ctabg \over J_{g-k}},
 \label{eqn:Tgk}
 \end{equation}
so that $HFF^*_g=\bigoplus \L^k_0 \ox T_{g,k}$. Note that the ideals
$J_{g-k}$ depend (in principle) on $g$ and $k$, and not only on the
difference $g-k$. 
The discussion following proposition 5.1 in~\cite{hff} shows that
$T_{g,k}$ is generated, as a free $\Ct$-module, by $\a^a\b^b\g^c$, 
$a+b+c < g-k$. 

\begin{prop}
\label{prop:Tgk}
  There is a direct sum decomposition
 \begin{equation}
  T_{g,k}=\bigoplus_{r=-(g-k-1)}^{g-k-1} R_{g,k,r},
 \label{eqn:Rgkr}
 \end{equation}
  where $R_{g,k,r}$ are free $\Ct$-modules. For $r$ even, 
  the eigenvalues of $(\a,\b,\g)$ in  $R_{g,k,r}$ are of the form
  $(4r\ima +O(t),8+O(t),0)$.  For $r$ odd, 
  the eigenvalues of $(\a,\b,\g)$ in  $R_{g,k,r}$ are of the form
  $(4r +O(t),-8+O(t),0)$. Here $O(t)$ means any series $f(t)\in
  t\Ct$. 
\end{prop}

\begin{pf}
From~\cite[lemma 5.2]{hff} we have that 
$\g J_r \subset J_{r+1}$, for all $0 \leq r\leq g-1$. 
It follows that $\g^{g-k} \in J_{g-k}$, i.e.\ $\g^{g-k}=0$ in
$T_{g,k}$. 
So the only eigenvalue of $\g$ is zero. To find the eigenvalues of
$\a$ 
and $\b$ we restrict to study the quotient $\rT_{g,k}=T_{g,k}/\g
T_{g,k}$, 
which is the space $\overline\cF_{g-k}$ of~\cite[proposition
5.11]{hff}. 
By~\cite[lemma 5.12]{hff} there is an exact sequence
\begin{equation}
  \bigoplus_{-(g-k) \leq r\leq g-k \atop r\equiv g-k \pmod 2}
R_{g-k+1,r}
  \inc \rT_{g,k-1} \surj \rT_{g,k} 
\label{eqn:cFr}
\end{equation}
where $R_{g-k+1,r}$ is a free $\Ct$-module of rank $1$, such that for
$r$ 
even, $\a=4r\ima+O(t)$ and $\b=8+O(t)$, and for $r$ odd,
$\a=4r+O(t)$ and $\b=-8+O(t)$. Starting with $\rT_{g,g}=0$, we find
by descending induction that the eigenvalues of $\rT_{g,k}$ are as
in the statement. 

Now let $R_{g,k,r} \subset T_{g,k}$ be the $\Ct$-submodule
generated by all vectors $v\in T_{g,k}$ 
annihilated by some power of either $\a-(4r+f(t))$
if $r$ is odd, or $\a-(4r\ima+f(t))$ if $r$ is even, where $f(t)\in
t\Ct$
is a series (depending in principle on $v$). Then an application of
the
Chinese Remainder Theorem yields that
$T_{g,k}=\bigoplus\limits_{r=-(g-k-1)}^{g-k-1}
R_{g,k,r}$ over $\Ct$ (we only need to use that two polynomials 
$p_{1,t}(X),p_{2,t}(X) \in \Ct [X]$ are coprime over $\Ct$ if 
$p_{1,0}(X),p_{2,0}(X) \in \CC[X]$ have no common roots).
There is an alternative way of looking at this. Take the maximal
ideal
$\frm_r \subset \Ctabg$ given as $\frm_r=(t,\a-4r,\b-8,\g)$ if $r$ is
odd, 
$\frm_r=(t,\a-4r\ima,\b+8,\g)$ if $r$ is even. Then 
$R_{g,k,r}=(T_{g,k})_{\frm_r}$ is the localisation of $T_{g,k}$ at
$\frm_r$.
\end{pf}

Now we shall filter $T_{g,k}$ (and also each of the $R_{g,k,r}$) by
the ideals
generated by the powers of $\g$ and consider the associated graded
rings
 \begin{equation}
   \Gr_{\g} T_{g,k}= \bigoplus_{i \geq 0} {\g^i T_{g,k}\over
\g^{i+1}T_{g,k}}=
   \bigoplus_{r=-(g-k-1)}^{g-k-1} \Gr_{\g} R_{g,k,r}.
 \label{eqn:Gr}
 \end{equation}
The reason for doing this is that the associated graded rings are
always 
easier to describe than the rings themselves. 

\begin{lem}
\label{lem:uff}
  For any $i=0,1,\ldots, g-k-1$, there is a well-defined map 
  $$
  T_{g,k+i}/\g T_{g,k+i} \stackrel{\g^i}{\ar} \g^i T_{g,k}/\g^{i+1}
  T_{g,k},
  $$
  which moreover is an isomorphism.
\end{lem}

\begin{pf}
  From~\cite[lemma 5.2]{hff} we have $\g J_r \subset J_{r+1} \subset
  J_r$, for all $1 \leq r\leq g-1$, so that $\g^i J_{g-k-i} \subset
  J_{g-k}$ and
  the map of the statement is well-defined. Surjectivity is obvious.
To prove
  injectivity we work as follows. Note that $\a^a\b^b$, $a+b <g-k-i$,
are a 
  basis for $T_{g,k+i}/\g T_{g,k+i}$. Then $\a^a\b^b\g^i$, $a+b
<g-k-i$, 
  generate the $\Ct$-module $\g^i T_{g,k}/\g^{i+1} T_{g,k}$. So
$\rk_{\Ct}
  (\g^i T_{g,k}/\g^{i+1} T_{g,k}) \leq {g-k-i+1 \choose 2}$ and 
$$
  \rk_{\Ct} T_{g,k} = \rk_{\Ct} (\Gr_{\g} T_{g,k}) = \sum \rk_{\Ct}
\left(
  {\g^i T_{g,k} \over \g^{i+1} T_{g,k}} \right)
  \leq \sum {g-k-i+1 \choose 2} =\rk_{\Ct} T_{g,k}.
$$
  Therefore equality must hold, $\rk_{\Ct} (\g^i T_{g,k}/\g^{i+1}
T_{g,k}) = 
  {g-k-i+1 \choose 2}$, and $\a^a\b^b\g^i$, $a+b <g-k-i$, are a basis
for
  $\g^i T_{g,k}/\g^{i+1} T_{g,k}$. This completes the proof.
\end{pf}

Lemma~\ref{lem:uff} gives that, as $\Ctab$-algebras,
\begin{equation}
  \Gr_{\g} T_{g,k} \iso \bigoplus_{i= 0}^{g-k-1} \g^i \cdot
{T_{g,k+i}
  \over \g T_{g,k+i}}.
\label{eqn:Grg}
\end{equation}
Thus in order to describe~\eqref{eqn:Grg} we only need to understand
the $\Ctab$-module
\begin{equation}
\rT_{g,k}= T_{g,k}/\g T_{g,k},
\label{eqn:rTgk}
\end{equation}
for $0\leq k\leq g-1$. In~\cite{hff} we have a description
of $\rT_{g,k}$ (where it is denoted by $\overline{\cF}_{g-k}$).

\begin{prop} {\rm (\cite[proposition 5.11]{hff})}
\label{prop:rT}
  For all $0\leq k\leq g-1$, $\rT_{g,k}$ is a free $\Ct$-module
  of rank $\rk_{\Ct} \rT_{g,k}={g-k+1\choose 2}$, 
  with basis $\a^a\b^b$, $a+b<g-k$. Moreover
  $\rT_{g,k}=\Ctab/\rJ_{g-k}$, with $\rJ_r=(\rR^1_r, \rR^2_r)$, where
  $\rR^1_0=1$, $\rR^2_0=0$, and for all $0 \leq r \leq g-1$,
  $$
   \left\{ \begin{array}{l} \rR_{r+1}^1 = (\a+\bar f_{12}) 
   \rR_r^1+ r^2(1+\bar f_{22})\rR_r^2 
   \\ \rR_{r+1}^2 = (\b+(-1)^{r+1}8+\bar f_{21}) \rR_r^1 +\bar f_{22}
\rR_r^2 
   \end{array} \right. 
  $$ 
  for some (unknown) $\bar f_{ij} \in t\Ctab$, dependent on $r$ and
$g$.
  \hfill $\Box$
\end{prop}

By proposition~\ref{prop:Tgk} and equation~\eqref{eqn:rTgk}
there is a decomposition
\begin{equation}
 \rT_{g,k} = \bigoplus_{r=-(g-k-1)}^{g-k-1} \rR_{g,k,r}, \qquad
  \text{where } \rR_{g,k,r}={R_{g,k,r} \over \g R_{g,k,r}}.
\label{eqn:rRgkr}
\end{equation}
Here $\rR_{g,k,r}$ is an algebra over $\Ct$, 
free of finite rank as $\Ct$-module. Indeed
$\rR_{g,k,r}$ is characterized as the subset of $\rT_{g,k}$ where
the eigenvalues of $(\a,\b)$ are of the form $(4 r + O(t),-8+O(t))$
if $r$ is odd, $(4r\ima +O(t),8+O(t))$ if $r$ is even. 
In $\rR_{g,k,r}$, we shall put $\rb=\b +(-1)^{r+1} 8$. 
We may consider $\rR_{g,k,r}$ only as a $\Ctrb$-algebra. Then 
we have the following result

\begin{lem}
\label{lem:fin}
  The rank of $\rR_{g,k,r}$ is $d+1=\left[{g-k-1- |r| \over 2} 
  \right] +1$, for $-(g-k-1)\leq r \leq g-k-1$. That is,
  $\rR_{g,k,r}=\Ctrb/(P_{d+1,t}(\rb))$, where $P_{d+1,t}(\rb) \in
  \Ctrb$ is a monic
  polynomial of degree $d+1$, dependent on $g$ and $r$, 
  such that $P_{d+1,0}=\rb^{d+1}$.
\end{lem}

\begin{pf}
Decomposing~\eqref{eqn:cFr} according to the eigenvalues of $\a$ we
get the exact sequence
\begin{equation}
  R_{g-k+1,r} \inc \rR_{g,k-1,r} \surj \rR_{g,k,r},
\label{eqn:cFr2}
\end{equation}
for $-(g-k-1)\leq r \leq g-k-1$,
where we put $R_{g-k+1,r}=0$ if $r\not\equiv g-k\pmod 2$. If 
$r \equiv g-k\pmod 2$ then $R_{g-k+1,r}$ is a free $\Ct$-module of
rank $1$ such that $\b=(-1)^r8+O(t)$. Therefore
$$
  \left\{ \begin{array}{lcll}
  \rk_{\Ct} \rR_{g,k-1,r} &= &\rk_{\Ct} \rR_{g,k,r} +1, \qquad & 
   r \equiv g-k \pmod 2 \\
  \rk_{\Ct} \rR_{g,k-1,r} &=& \rk_{\Ct} \rR_{g,k,r}, &r\not\equiv g-k
  \pmod 2
  \end{array} \right.
$$
Also~\eqref{eqn:cFr2} for $r=\pm (g-k)$ 
gives that $\rk_{\Ct} \rR_{g,g-|r|-1,r}=1$, so 
we get by descending induction on $k$ that 
$\rk_{\Ct} \rR_{g,k,r}=\left[{g-k-1-|r| \over 2}\right]+1$, 
for $0\leq k \leq g-|r|-1$. To prove the second assertion, 
suppose that we have already $\rR_{g,k,r}=\Ctrb/(P_{d+1,t}(\rb))$, 
where $P_{d+1,t}(\rb) \in \Ctrb$ is a monic
polynomial of degree $d+1$ such that $P_{d+1,0}(\rb)=\rb^{d+1}$,
and let us prove the result for $\rR_{g,k-1,r}$.
If $r \not\equiv g-k \pmod 2$ then $\rR_{g,k-1,r} \iso \rR_{g,k,r}=
\Ctrb/(P_{d+1,t}(\rb))$, so we are done. 
If $r \equiv g-k \pmod 2$ then let
$(-1)^r8+f(t)$, with $f(t)\in t\Ct$, be the eigenvalue of $\b$ 
in $R_{g-k+1,r}$, so that $R_{g-k+1,r}=\Ctrb/(\rb+f(t))$. 
If the exact sequence~\eqref{eqn:cFr2} 
is not split (over $\Ct$) then $\rR_{g,k-1,r}=
\Ctrb/(P_{d+2,t}(\rb))$, where $P_{d+2,t}(\rb)= 
(\rb+f(t))P_{d+1,t}(\rb)$, and the result follows.

Let us see that~\eqref{eqn:cFr2} does not split. Consider
as in~\cite[proposition 4.3]{hff}, the spaces $F_r=(HF^*_r)_I$
and $\overline{F}_r=F_r/\g F_r$, so that
the map given by equating $t=0$ is (see~\cite[proposition 5.11]{hff}), 
\begin{equation}
  \rT_{g,k-1}=\overline{\cF}_{g-k+1} \surj \overline{F}_{g-k+1}.
\label{eqn:q}
\end{equation}
By~\cite[proposition 20]{hf}, there is a decomposition 
$\overline{F}_{g-k+1} = \bigoplus\limits_{-(g-k)\leq i\leq g-k} A_i$, where 
$\a- 4i\ima$ and $\b-8$ are nilpotent for $i$ even, 
$\a- 4i$ and $b+8$ are nilpotent for $i$ odd. Also from the proof 
of~\cite[proposition 20]{hf} we have that
$$
  \left\{ \begin{array}{l}   
  \overline{F}_g/(\b+8)=\CC[\a]/\left((\a^2 +(2[{g-1\over 2}])^2 16) 
  \cdots (\a^2 + 2^2 16)\a \right)  \\
  \overline{F}_g/(\b-8)=\CC[\a]/\left((\a^2 -(2[{g\over 2}]-1)^2 16) 
  \cdots (\a^2 - 1^2 16) \right)
  \end{array} \right.
  $$
So in $A_i$, $\a- 4i\ima$ is a polynomial in $\b-8$ (with no independent 
term) for $i$ even, $\a- 4i$ is a polynomial in $\b+8$ (with no independent 
term) for $i$ odd. In particular $A_r=\Crb/(\rb^{d+2})$ as $\Crb$-algebra.
By~\eqref{eqn:q}, we have $\rR_{g,k-1,r} \surj A_r$ and it follows
that $\rR_{g,k-1,r}$ cannot be a direct sum.
\end{pf}

We summarize what we have in the following

\begin{thm}
\label{thm:main}
  The Fukaya-Floer homology of $Y=\Y$, where $\S$ is a surface of
  genus $g\geq 1$, is
  $$
    HFF^*_g = \bigoplus_{0\leq k \leq g-1 \atop -(g-k-1) \leq r \leq
    g-k-1} \L^k_0(\seq{\q}{1}{2g}) \otimes R_{g,k,r},
  $$
  where $R_{g,k,r}$ is a $\Ctabg$-algebra, free as $\Ct$-module,
characterized
  by the condition that the eigenvalues of $\a$ in $R_{g,k,r}$
  are $4r+O(t)$ if $r$ is odd, $4r\ima+O(t)$ if $r$ is even.
  In $R_{g,k,r}$ put $\rb=\b+(-1)^{r+1}8$. Then as $\Ctrb$-module, 
$$
  \Gr_{\g} (R_{g,k,r}) \iso \bigoplus_{i=0}^{g-k-|r|-1} \g^i \cdot 
  {\Ctrb \over  (P_{d(g,r,k,i)+1,t}(\rb))},
$$
  where $d=d(g,r,k,i)=[{g-k-|r|-i-1 \over 2}]$ and where
  $P_{d+1,t}(\rb) \in \Ctrb$ is a monic
  polynomial of degree $d+1$ dependent on $g$ and $r$ such
  that $P_{d+1,0}(\rb)=\rb^{d+1}$.\hfill $\Box$ 
\end{thm}

We may decompose the Fukaya-Floer homology with respect to the
eigenvalues of
$\a$ alone,
\begin{equation}
  HFF^*_g =\bigoplus\limits_{r=-(g-1)}^{g-1} \cH_r,
\label{eqn:Hr}
\end{equation}
where the eigenvalues of $\a$ in $\cH_r$ are of the form $4r+O(t)$ if 
$r$ is odd,
$4r\ima+O(t)$ if $r$ is even. Note that~\eqref{eqn:Hr} is actually an
orthogonal decomposition. Theorem~\ref{thm:main} implies that
\begin{equation}
  \cH_r= \bigoplus_{k=0}^{g-|r|-1} \L^k_0 \ox R_{g,k,r},
\label{eqn:Hr2}
\end{equation}
which can be considered as a $\Ctrbg$-module. 
So there is a $\Spz$-equivariant epimorphism
\begin{equation}
  \tilde{\AA}(\S) \otimes \Ct \iso \L^* (\seq{\q}{1}{2g}) \ox \Ctrb 
  \surj \cH_r=\bigoplus_{k=0}^{g-|r|-1} \L^k_0 \ox R_{g,k,r}\, ,
\label{eqn:impo}
\end{equation}
where the isomorphism of the left hand side is given by
$\g_i \mapsto \q_i$, $\wp \mapsto \rb$, and endows 
$\L^* \ox \Ctrb$
with a grading $d$ with $d(\rb)=2$ and $d(\q_i)=1$, $1\leq i \leq 2g$
(the elements of $\Ct$ are the coefficients, with degree $0$). 
$\cH_r$ is not graded in general, as the kernel of~\eqref{eqn:impo}
is not a graded ideal.

\begin{cor}
\label{cor:relat}
  For $-(g-1) \leq r \leq g-1$ we have a presentation
  $$ 
  \cH_r =\bigoplus_{k=0}^{g-|r|-1} \L^k_0 \ox {\Ctrbg \over \cI_k},
  $$
  where $\rb=\b+(-1)^r8$. A basis for $\Ctrbg /\cI_k$ is given by
  $\rb^i\g^j$, $2i+j <g-k-|r|$. For $0\leq k\leq g-|r|$, there are 
  polynomials (depending on $g$ and $r$)
  $$
   R_{k}=  P_{d+1,t}(\rb) - \sum_{2i+j < g-k-|r| \atop j > 0}
    c^{k}_{ij}\rb^i\g^j,
  $$
  where $d=[{g-k-|r|-1 \over 2}]$, $P_{d+1,t}(\rb)\in \Ctrb$ is a
  monic polynomial of degree $d+1$ such that $P_{d+1,0}(\rb)=\rb^{d+1}$,
  and $c^{k}_{ij} \in \Ct$, satisfying
  $\cI_k=(R_k ,\g R_{k+1},\g^2 R_{k+2},\ldots, \g^{g-k-|r|})$.
  (Note that $R_{g-|r|}= 1$.)
\end{cor}

\begin{pf}
  Theorem~\ref{thm:main} implies that the elements $\g^i\rb^j$, with
  $j \leq d(g,r,k,i)=[{g-k-|r| -i-1 \over 2}]$ form a basis for
$R_{g,k,r}$.
  This condition is equivalent to $2j+i < g-k-|r|$.
  The ideal $\cI_k$ is generated by polynomials
  \begin{equation}
   P_{n,t}(\rb) \g^m - \sum_{2i+j < g-k-|r| \atop j > m} c^{nm}_{ij} 
  \rb^i\g^j,
  \label{eqn:esta}
  \end{equation}
  for some $c^{nm}_{ij} \in \Ct$, where $0\leq m \leq g-k-|r|$ and 
  $n=[{g-k-|r|-m-1 \over 2}]$+1. 
  Now let $R_k$ be the relation~\eqref{eqn:esta} with $m=0$.
  The inclusions $\g \cI_{k+1} \subset \cI_{k} \subset \cI_{k+1}$
yield 
  that the relation~\eqref{eqn:esta} 
  in $\cI_k$ equals $\g^m R_{k+m}$. The result follows.
\end{pf}

\section{Floer homology and symmetric products of $\S$}
\label{sec:hf}

Let $HF_g^*=HF^*(Y)$ be the Floer homology of $Y=\Y$, where $\S$ is a
surface
of genus $g\geq 1$, and for the $SO(3)$-bundle with $w_2=\PD [\SS^1]
\in H^2(Y;\ZZ/2\ZZ)$. This has
been determined by the works of Dostoglou and Salamon~\cite{DS} and
the ring 
structure by the author in~\cite{hf}. The Floer homology $HF_g^*$ is
a 
finite dimensional algebra over $\CC$ (with a graduation modulo $4$)
and
there is a natural epimorphism $HFF_g^* \surj HF_g^*$ given
by equating $t=0$ (see~\cite{hff}). As a consequence of
theorem~\ref{thm:main}
we have 

\begin{prop}
\label{prop:m}
  The Floer homology of $Y=\Y$, where $\S$ is a surface of genus 
  $g\geq 1$, is
  $$
    HF^*_g = \bigoplus_{0\leq k \leq g-1 \atop -(g-k-1) \leq r \leq
g-k-1}
    \L^k_0(\seq{\q}{1}{2g}) \otimes \hat{R}_{g,k,r},
  $$
  where $\hat{R}_{g,k,r}$ is a $\Cabg$-algebra, such that
  the eigenvalue of $\a$ in $\hat{R}_{g,k,r}$ is $4r$ if $r$ is 
  odd, $4r\ima$ if $r$ is even.
  In $\hat{R}_{g,k,r}$ put $\rb=\b+(-1)^{r+1}8$. Then as
$\CC[\rb]$-module, 
$$
  \Gr_{\g} (\hat{R}_{g,k,r}) \iso \bigoplus_{i=0}^{g-k-|r|-1} \g^i
\cdot 
  {\Ctrb \over (\rb^{d(g,r,k,i)+1}))},
$$
  where $d(g,r,k,i)=[{g-k-|r|-i-1 \over 2}]$.\hfill $\Box$ 
\end{prop}

  We have an artinian decomposition $HF^*_g
=\bigoplus\limits_{r=-(g-1)}^{g-1}
  H_r$, where the eigenvalue of $\a$ in $H_r$ is $4r$ if $r$ is odd, 
  $4r\ima$ if $r$ is even.
  For $-(g-1) \leq r \leq g-1$ there is a $\Spz$-equivariant
epimorphism
  $$
  \tilde{\AA}(\S) \iso \L^* (\seq{\q}{1}{2g}) \ox \CC[\rb] \surj H_r,
  $$
  where $\rb=\b+(-1)^r8$. The isomorphism of the left hand side is
given by
  $\g_i \mapsto \q_i$, $\wp \mapsto \rb$, and endows 
$\L^* \ox \Crb$
with a grading $d$ with $d(\rb)=2$ and $d(\q_i)=1$, $1\leq i \leq
2g$.
Let us see the interesting fact that the product of two elements of 
$H_r$ is a sum of elements of the same or higher degree.

\begin{lem}
\label{lem:d(e)}
Consider the local ring  $\hat{R}_{g,k,r}$ and put $\rg =g-k-|r|$. 
Then the elements $\rb^n\g^m$, $2n+m < \rg$, form a basis of 
$\hat{R}_{g,k,r}$. Any  $w=\rb^n\g^m$ with $2n+m \geq \rg$ 
can be written as
\begin{equation}
  w= \sum_{2i+j < \rg, j\geq m \atop i+j \geq n+m} c_{ij} \rb^i\g^j,
\label{eqn:e}
\end{equation}
for some $c_{ij} \in \CC$,
i.e.\ $w$ is a linear combination of monomials of the basis of the 
same or higher degree.
\end{lem}

\begin{pf}
  Let $w=\rb^n\g^m$ with $2n+m \geq \rg$. We shall prove that $w$
  can be
  written as~\eqref{eqn:e}, by descending induction on $k$. 
  For $k=g-|r|-1$, it is $\rg=1$, $\hat R_{g,k,r}=\Crbg/(\rb,\g)$ 
  and the statement is true. Now let $0 \leq k <g-|r|-1$ and suppose 
  that the statement is proved for $k+1$.
  Note that the inclusion $\g J_{g-k-1}\subset J_{g-k}$ yields a
  well defined map $T_{g,k+1} \stackrel{\g}{\ar} T_{g,k}$, which in
  turn gives maps  $R_{g,k+1,r} \stackrel{\g}{\ar} R_{g,k,r}$
  and $\hat{R}_{g,k+1,r} \stackrel{\g}{\ar} \hat{R}_{g,k,r}$.
  So if $m>0$, the inductive hypothesis implies that  
$$
  \rb^n\g^{m-1}= \sum_{2i+j < \rg-1, j\geq m-1 \atop i+j \geq n+m-1}
  a_{ij}\rb^i\g^j
$$
  in $\hat{R}_{g,k+1,r}$. Now multiplying by $\g$ we have the
  equation~\eqref{eqn:e}
  in $\hat{R}_{g,k,r}$ with $c_{ij}=a_{i,j-1}$ for $j\geq 1$,
  $c_{i0}=0$.

  For the case $m=0$, we see that it is enough to prove the statement
  for $\rb^{d+1}$, where  $d=[{\rg-1 \over 2}]$. This is so since
  once we have equation~\eqref{eqn:e} for $w=\rb^{d+1}$, we may 
  multiply by $\rb$ and use recurrence to get equation~\eqref{eqn:e} 
  for any $w=\rb^n$ with $n >d+1$.

  Now let $w=\rb^{d+1}$, and write it in terms of the basis as
  \begin{equation}
    \rb^{d+1}= \sum_{2i+j < \rg \atop j \geq 1 } c_{ij} \rb^i\g^j.
  \label{eqn:mia}
  \end{equation}
  Note that it must be $j\geq 1$ from proposition~\ref{prop:m}.
  Multiplying both sides of~\eqref{eqn:mia} by $\g$ we get
  $\rb^{d+1}\g= \sum c_{ij} \rb^i\g^{j+1}$. 
  Using the case already proved above (for $m>0$), we get that
 $$
  \sum_{2i+j < \rg, j\geq 1 \atop i+j+1 \leq d+1} c_{ij}
\rb^i\g^{j+1}
 $$
  is expressible as a linear combination of monomials of the basis 
  of degree bigger
  or equal than $d+2$. This is impossible since $\rb^i\g^{j+1}$, 
  $i+j+1 \leq d+1$, $j \geq 1$, are themselves monomials of the
basis. Indeed,
  $2i+(j+1) \leq 2d+1-j \leq 2d =2[{\rg-1 \over 2}] <\rg$.
  So it must be $c_{ij}=0$ in~\eqref{eqn:mia} whenever $i+j \leq d$.
\end{pf}

Analogously to corollary~\ref{cor:relat} we have in the case
of the Floer homology $HF^*_g$ the following result

\begin{prop}
\label{prop:hf}
  We have an artinian decomposition $HF^*_g
  =\bigoplus\limits_{r=-(g-1)}^{g-1}
  H_r$, where the eigenvalue of $\a$ in $H_r$ is $4r$ if $r$ is odd, 
  $4r\ima$ if $r$ is even.
  For $-(g-1) \leq r \leq g-1$ there is a $\Spz$-equivariant
  epimorphism
  $\tilde{\AA}(\S) \iso \L^* (\seq{\q}{1}{2g}) \ox \CC[\rb] \surj H_r$,
  where $\rb=\b+(-1)^r8$. There is a presentation 
  $$ 
  H_r =\bigoplus_{k=0}^{g-|r|-1} \L^k_0 \ox {\CC[\rb,\g] \over I_k}. 
  $$
  Also for $0 \leq k\leq g-|r|$, there are polynomials 
  $$
   R_{k}= \rb^{d+1} - \sum_{2i+j < g-k-|r| \atop i+j \geq d+1,\,j>0} 
  c^{k}_{ij} \rb^i\g^j,
  $$
  where $d=[{g-k-|r|-1 \over 2}]$ and $c^{k}_{ij} \in \CC$. 
  Then $I_k=(R_k ,\g R_{k+1},\g^2 R_{k+2},\ldots, \g^{g-k-|r|})$. 
\hfill $\Box$
\end{prop}

\begin{rem}
  In proposition~\ref{prop:hf} we have not taken any effort in 
  looking for a minimal set of generators of the ideals $I_k$. 
  We conjecture that the first two generators suffice. 
\end{rem}

It is our intention to relate the pieces $H_r$ with the cohomology of
the 
symmetric product $s^{g-|r|-1}\S$ of the surface $\S$, as predicted 
in~\cite[conjecture 24]{hf}.
Fix $0 \leq d\leq g-1$. We shall review the description of
$H^*(s^d\S)$ 
from~\cite{macdonald}, putting it into the right form for our
purposes.
We interpret $s^d\S$ as the moduli space of degree $d$ effective
divisors on
$\S$ (for this we need to put a complex structure on $\S$, but at the
end 
this will be irrelevant for the description of the cohomology). 
Let $D \subset s^d\S \x\S$ be the universal divisor. Then we define
$$
  \left\{ \begin{array}{l}
  \q_i =c_1(D) /\g_i  \in H^1(s^d\S), \qquad  1\leq i \leq 2g \\
  \eta = c_1(D) / x \in H^2(s^d\S) \end{array} \right.
$$
These elements generate $H^*(s^d\S)$. More precisely, there is a
(graded!)
$\Spz$-equivariant epimorphism
\begin{equation}
  \tilde{\AA}(\S) \iso \L^*(\seq{\q}{1}{2g})\otimes\CC[\eta] \surj
H^*(s^d\S).
\label{eqn:00}
\end{equation}
In particular (for $d>0$) $H^1(s^d\S)$ is freely generated by 
$\seq{\q}{1}{2g}$. We
put $\L^k_0=\L^k_0(\seq{\q}{1}{2g})$. Note that since $H^*(s^d\S)$ is
graded 
from $0$ to $2d$, $\L^k_0$ goes to zero under~\eqref{eqn:00} for
$k>d$.
Finally note that $\eta$ and $\h=\sum \q_i\q_{g+i}$ generate the
invariant
part $H^*(s^d\S)_I$. 

\begin{prop}
\label{prop:hsdS}
  For $-(g-1) \leq r \leq g-1$ there is a presentation 
  $$ 
    H^*(s^d\S) =\bigoplus_{k=0}^{d} \L^k_0 \ox {\CC[\eta,\h] \over
J_k}, 
  $$
  with $J_k=(R_k ,\h R_{k+1},\h^2 R_{k+2},\ldots, \h^{d+1-k})$ and
  $$
  R_k= \sum\limits_{i=0}^{\a} {{d-k-\a+1 \choose i} 
   \over{g-k \choose i}} {(-\h)^i \over i!}\eta^{\a-i},
  $$
  where $\a=[{d-k \over 2}]+1$, for $0 \leq k \leq d$, and
$R_{d+1}=1$. 
\end{prop}

\begin{pf}
  The relations of $H^*(s^d\S)$ are given in~\cite{macdonald} and are
  the following
  \begin{equation}
    \eta^r \prod_{i\in I} (\eta - \q_i \q_{g+i}) \prod_{j \in J} \q_j 
    \prod_{k\in K} \q_{g+k},
  \label{eqn:j2}
  \end{equation}
  for $I, J, K \subset\{1,\ldots, g\}$ disjoint and $r+2|I|+|J|+|K|
\geq d+1$. 
  Suppose $d \not{\!\!\equiv\>} k \pmod 2$, so $d+1=2\a+k$. Take the
relation 
  $\q_1 \cdots \q_k \prod_{i\in I} (\eta - \q_i \q_{g+i})$, with 
  $I\subset \{ k+1,\ldots, g\}$, $|I|=\a$. Let
$\text{Sp}\,(2g-2k,\ZZ)$ act 
  on the relation, acting in the standard way on $\{\q_{k+1} ,\ldots,
\q_g, 
  \q_{g+k+1},\ldots, \q_{2g}\}$. This produces the relation
  $$
   \q_1 \cdots \q_k   \sum_{i=0}^{\a} {{\a \choose i} \over{g-k
\choose i} }
  \eta^{\a-i} {(-\h')^i \over i!},
  $$
  where $\h'=\sum_{i=k+1}^g \q_i\q_{g+i}$. Using $\a=d-k-\a+1$ 
  this is equivalent to 
  the relation $\q_1 \cdots \q_k R_k$, i.e. $R_k \in J_k$. 
  For $d \equiv k \pmod 2$ take~\eqref{eqn:j2} with $r=1$,
$J=\{1,\ldots, k\}$,
  $K=\emptyset$, $d+1=2(\a-1)+k+1$, and proceed as above to get 
  $$
   \q_1 \cdots \q_k   \sum_{i=0}^{\a} {{\a-1 \choose i} \over{g-k
\choose i} }
  \eta^{\a-1-i} {(-\h)^i \over i!}.
  $$
  Then we use $\a-1=d-k-\a+1$ to get $R_k \in J_k$.
 
  From the inclusions $\h J_{k+1} \subset J_k \subset J_{k+1}$ we get
  that
  $\h R_{k+1}, \ldots, \h^{d-k} R_d,\h^{d+1-k} \in J_k$. So 
  $(\h R_{k+1}, \ldots, \h^{d-k} R_d,\h^{d+1-k}) \subset J_k$ and 
  there is an epimorphism
  \begin{equation}
    A^*={\CC[\eta,\h] \over (R_k, \h R_{k+1},\ldots, \h^{d+1-k})}
    \surj B^*={\CC[\eta,\h] \over J_k}.
  \label{eqn:j3}
  \end{equation}
  Both $A^*$ and $B^*$ are (evenly) graded rings and~\eqref{eqn:j3}
is graded.
  As $H^*(s^d\S)$ has Poincar\'e duality with $\dim s^d\S=2d$, we
conclude that
  $B^*$ has Poincar\'e duality $B^* \otimes B^{2d-2k-*} \ar \CC$.
Also all relations
  in~\eqref{eqn:j2} have degrees bigger than or equal to $d+1$, so
all elements in
  $J_k$ have degrees bigger than or equal to $d-k+1$. The generators
of the ideal of
  relations of $A^*$ have also degrees bigger than or equal to
$d-k+1$, so
  $\dim A^{2i}=\dim B^{2i}$, for $0 \leq 2i \leq d-k$. Now 
  $$ 
  \Gr_{\h} A^*=\bigoplus_{i=0}^{d-k} \h^i \cdot {\CC[\eta] \over
(\eta^{[(d-k-i)/2]+1}) }.
  $$
  Hence for $d-k \leq 2i \leq 2d-2k$, we have $\dim A^{2i} = d-k-i+1
=\dim A^{2d-2k-2i}
  =\dim B^{2d-2k-2i}=\dim B^{2i}$. Therefore~\eqref{eqn:j3} is an
isomorphism.
\end{pf}

\begin{rem}
\label{rem:j4}
  In proposition~\ref{prop:hsdS}, the first two generators of $J_k$
suffice, i.e.\
  $J_k=(R_k, \h R_{k+1})$. This follows from
$$
  \left\{ \begin{array}{ll}
  R_k =R_{k+1} - {g-k-\a \over (g-k)(g-k-1)} \h R_{k+2}, &
d\not{\!\!\equiv\>} k \pmod 2 \\
  R_k =\eta R_{k+1} + {\a-1 \over (g-k)(g-k-1)} \h R_{k+2},\qquad & 
  d \equiv k \pmod 2 \end{array} \right.
$$
  where $\a=[{d-k \over 2}]+1$.
\end{rem}

As a consequence of proposition~\ref{prop:hf} we have
  $$
  \Gr_{\g} H_r = \bigoplus_{k=0}^{g-|r|-1} \L^k_0 \ox
  \Big( \bigoplus_{i=0}^{g-|r|-1-k} \g^i \cdot {\CC[\rb] \over 
  (\rb^{d(g,r,k,i)+1})} \Big),
  $$
where $d(g,r,k,i)=[{g-k-|r|-i-1 \over 2}]$. 
And as a consequence of proposition~\ref{prop:hsdS},
  $$ 
  \Gr_{\h} H^*(s^d\S) =\bigoplus_{k=0}^{d} \L^k_0 \ox
  \Big( \bigoplus_{i=0}^{d-k} \h^i \cdot {\CC[\eta] \over 
  (\eta^{[(d-k-i)/2]+1}) } \Big).
  $$
This proves the following

\begin{cor}
\label{cor:hom-symm}
  There is an isomorphism
  $\Gr_{\g} H_r \iso \Gr_{\h} H^*(s^{g-|r|-1}\S)$, for $-(g-1) \leq
  r\leq g-1$. \hfill $\Box$
\end{cor}

In particular $H_r$ and $H^*(s^{g-|r|-1}\S)$ are isomorphic as 
$\Spz$-representations, which proves~\cite[conjecture 24]{hf}. 
We believe that in fact $H_r$ is isomorphic to the quantum cohomology 
$QH^*(s^{g-|r|-1}\S)$ (see~\cite{BT} for a partial computation of the
latter).

\section{New relations for the Fukaya-Floer homology}
\label{sec:relations}

In this section we are going to get more information on the shape of
the relations for $\cH_r$ given in corollary~\ref{cor:relat} 
in order to prove our main theorems. Our next result is a weaker
version of lemma~\ref{lem:d(e)} for the Fukaya-Floer homology,
telling us that the degrees homogeneous components 
of the relations for $\cH_r$ cannot be too low.

\begin{lem}
\label{lem:e(e)}
  With the notations of corollary~\ref{cor:relat}, let $e$ be the 
  multiplicity of the root $\rb=0$ in $P_{d+1,t}(\rb)$. Then we have
$$
  R_{k}= P_{d+1,t}(\rb) - \sum_{2i+j<g-k-|r|\atop i+j \geq e, j > 0}
  c^{k}_{ij} \rb^i\g^j.
$$
\end{lem}

\begin{pf}
  Let $d=d(g,r,k,0)$ and $\rg=g-k-|r|$. Write 
  $P_{d+1,t}(\rb)=\rb^e Q_t(\rb)$, where
  $Q_t(0)\neq 0$. We shall prove by descending induction on $k$ the
  following two statements: 
\begin{enumerate}
  \item $P_{d+1,t}(\rb)= \sum\limits_{2i+j<\rg\atop i+j \geq e, j > 0} 
  c_{ij} \rb^i\g^j \in R_{g,k,r}$, for some $c_{ij}\in \Ct$.
  \item $R_{g,k,r}=\bigoplus\limits_{i+j<e} \la \rb^i\g^j \ra \bigoplus
  \la \rb^i\g^j | i+j \geq e\ra$, i.e.\ all non zero monomials 
  of all relations in $R_{g,k,r}$ have degrees bigger or equal 
  than $e$.
\end{enumerate}
 
  For $k=g-|r|-1$ the statement is obvious. 
  Now let $0 \leq k <g-|r|-1$ and suppose that the statement 
  is proved for any number strictly bigger than $k$. 
  By corollary~\ref{cor:relat} we may write $P_{d+1,t}(\rb)=
  \sum\limits_{2i+j<\rg, j > 0} 
  a_{ij} \rb^i\g^j$, in $R_{g,k,r}$. We have to prove that
  $a_{ij}=0$ for $i+j<e$.
  We have the following three cases:
\begin{itemize}
\item $d(g,r,k+1,0)=d$, i.e.\ $\rg$ even. Take the 
  relation $P_{d+1,t}(\rb)= \sum\limits_{2i+j<\rg-1 
  \atop i+j \geq e, j > 0} 
  c_{ij} \rb^i\g^j$ in $R_{g,k+1,r}$. The map $\g:R_{g,k+1,r}
  \inc R_{g,k,r}$ gives that 
  $$
  \g P_{d+1,t}(\rb)= \sum_{2i+j<\rg-1 \atop i+j \geq e, j > 0} 
  c_{ij} \rb^i\g^{j+1}= \sum_{2i+j<\rg, j > 0} a_{ij} \rb^i\g^{j+1},
  $$
  in $R_{g,k,r}$. Using that $\bigoplus\limits_{i+j<e} 
  \la \rb^i\g^{j+1} \ra \bigoplus \la \rb^i\g^{j+1} | 
  i+j \geq e\ra =\g R_{g,k+1,r} \subset R_{g,k,r}$,  
  we get that $a_{ij}=0$ for $i+j<e$, as required. 
  The second item also follows, for all the relations
  $\g^m R_{k+m}$, $m>0$, in $R_{g,k,r}$, involve terms of degrees
  bigger or equal than $e$ by induction hypothesis, and
  the new relation, $\rb^e Q_t(\rb)= \sum\limits_{i+j\geq e, j>0} a_{ij}
  \rb^i\g^j$, also involves terms of degrees bigger or
  equal than $e$.
\item $d(g,r,k+1,0)=d-1$ and $P_{d+1,t}(\rb)=(\rb+f(t)) P_{d,t}(\rb)$,
  with $f(t)\in t\Ct$ non zero. As above, the relation in
  $R_{g,k+1,r}$ produces 
  $\g P_{d,t}(\rb)= \sum\limits_{i+j \geq e, j > 0} 
  c_{ij} \rb^i\g^{j+1}$ in $R_{g,k,r}$. Multiplying by $\rb+f(t)$,
  $$
  \g P_{d+1,t}(\rb)= \sum_{i+j \geq e, j > 0} 
  c_{ij}' \rb^i\g^{j+1}= \sum_{2i+j<\rg, j > 0} 
  a_{ij} \rb^i\g^{j+1}
  $$ 
  in $R_{g,k,r}$. The rest is as above.
\item $d(g,r,k+1,0)=d-1$ and $P_{d+1,t}(\rb)=\rb P_{d,t}(\rb)$.
  Now the relation in $R_{g,k+1,r}$ is 
  $P_{d,t}(\rb)= \sum\limits_{i+j \geq e-1, j > 0} 
  c_{ij} \rb^i\g^j$, since the multiplicity of the root $\rb=0$ in 
  $P_{d,t}(\rb)$ is $e-1$. Multiplying by $\rb$,
 \begin{equation}
  \g P_{d+1,t}(\rb)= \sum_{i+j \geq e, j > 0} c_{ij}' 
  \rb^i\g^{j+1}= \sum_{2i+j<\rg, j > 0} a_{ij} \rb^i\g^{j+1},
 \label{eqn:lo}
 \end{equation}
  in $R_{g,k,r}$. This time we need to use the inclusion
  $\g^2: R_{g,k+2,r} \inc
  R_{g,k,r}$, which yields that $\bigoplus\limits_{i+j<e-1} 
  \la \rb^i\g^{j+2} \ra \bigoplus \la \rb^i\g^{j+2} 
  | i+j \geq e-1 \ra =\g^2 R_{g,k+2,r}
  \subset R_{g,k,r}$. Applying this to~\eqref{eqn:lo} we get that
  $a_{ij}=0$ for $i+j<e$. The second item follows as above.
\end{itemize}
\end{pf}

Now we aim to prove the vanishing of polynomials of high degree on
$\rb$ and $\g$, giving a explicit bound on such degree. 
If we want to have this for $\cH_r$ we must prove first 
that $P_{d+1,t}(\rb)=\rb^{d+1}$, this being the only case in which
$\rb$ is nilpotent in $\cH_r$. This is equivalent
to proving that the only eigenvalues of $\b$ on $HFF^*_g$ are $\pm 8$.
(We note also that if this is the case then the proof of
lemma~\ref{lem:d(e)} carries over to the Fukaya-Floer setting giving
a simpler method for proving lemma~\ref{lem:e(e)}.)
On the other hand, from the point of view of Donaldson invariants we
may well restrict attention to the effective Fukaya-Floer homology,
where the eigenvalues of $\b$ are already $\pm 8$ (see~\cite[theorem
5.13]{hff}). We recall its definition

\begin{defn} {\rm (\cite[definition 5.10]{hff})}
\label{def:f5.eHFF}
  We define the effective Fukaya-Floer homology as the 
  sub-$\Ct$-module
  $\whff \subset HFF^*_g$ generated by all $\p^w(X_1, z_1 e^{tD_1})$,
  for all $4$-manifolds $X_1$ with boundary $\bd X_1=Y=\Y$ such that
  $X=X_1\cup_Y A$ has $b^+>1$, $z_1 \in \AA(X_1)$, 
  $D_1 \subset X_1$ with $\bd D_1=\SS^1$
  and $w\in H^2(X_1;\ZZ)$ with $w|_Y=w_2=\PD[\SS^1]$.
\end{defn}

Clearly there is a decomposition
$\whff =\bigoplus\limits_{r=-(g-1)}^{g-1} \tilde{\cH}_r$,
where $\tilde{\cH}_r\subset \cH_r$. By~\cite[theorem 5.13]{hff}
the only eigenvalues of $(\rb,\g)$ on $\tilde{\cH}_r$ are $(0,0)$.
Decompose $\tilde{\cH}_r=\bigoplus\limits_{k=0}^{g-|r|-1} \L^k_0 \ox
\tilde{R}_{g,k,r}$, where $\tilde{R}_{g,k,r} \subset R_{g,k,r}$.
By lemma~\ref{lem:e(e)} and corollary~\ref{cor:relat}, 
\begin{equation}
  \g^m P_{d+1,t}(\rb)=\g^m \rb^e Q_t(\rb)=
  \sum_{2i+j<g-k-|r|-m \atop i+j\geq e, j>0} c_{ij} \rb^i\g^{j+m}
\label{eqn:ee}
\end{equation}
in $R_{g,k,r}$, where $d=[{g-k-|r|-m-1 \over 2}]$, $e$ is
the multiplicity of the root $\rb=0$ in $P_{d+1,t}(\rb)$
and $Q_t(\rb)$ does not have $\rb=0$ as root.
On $\tilde{R}_{g,k,r}$ the only eigenvalue of $\rb$ is zero, 
so $Q_t(\rb)$ is invertible and hence~\eqref{eqn:ee} produces
$$
  \g^m \rb^e =\sum_{i+j\geq e, j>0} c_{ij} \rb^i\g^{j+m},
$$
for some (different from the previous ones) coefficients
$c_{ij} \in \Ct$, as endomorphisms
acting on $\tilde{R}_{g,k,r}$. As $e\leq d+1$, this implies that
\begin{equation}
  \g^m \rb^{d+1} =\sum_{i+j\geq d+1, j>0} c_{ij} \rb^i\g^{j+m},
\label{eqn:efe}
\end{equation}
acting on $\tilde{R}_{g,k,r}$. 

\begin{cor}
\label{cor:e}
  Let $b=\rb^n \g_{i_1} \cdots \g_{i_m} \in \tilde{\AA}(\S)$. 
  If $d(b)=2n+m > 2(g-|r|-1)$ then $b=0$ acting on 
  $\tilde{\cH}_r$.
\end{cor}

\begin{pf}
  Let us first prove that any $b=\rb^n\g^m$ with $n+m\geq g-k-|r|$ 
  is zero acting on $\tilde{R}_{g,k,r}$. As $\g^{g-k-|r|}=0$ in 
  $R_{g,k,r}$ we may suppose $m<g-k-|r|$ and $n>0$.
  Clearly $n \geq d+1=[{g-k-|r|-m-1 \over 2}]+1$, 
  so we use~\eqref{eqn:efe} 
  to write the endomorphism $\rb^{n}\g^m$ as a
  $\Ct$-linear combination of endomorphisms $\rb^i\g^j$ with 
  $i+j\geq n+m$ and $j>m$. By recursion the claim follows.

  Now let $b=\rb^n\g_{i_1}\cdots\g_{i_m}$ with $2n+m>2g-2|r|-2$. 
  We have the simple fact
$$
  \g_{i_1} \cdots \g_{i_m} \L^k_0 \subset \L^{k+m}_0 \oplus
  \g\L^{k+m-2}_0 \cdots \oplus \g^m \L^{k-m}_0,
$$
  where a negative exponent is understood as that the corresponding
  term does not appear. Hence to see that $b$ kills
  $\L^0_k \tilde{R}_{g,k,r}$ it is enough to prove that $\rb^n\g^i=0$ on 
  $\tilde{R}_{g,k+m-2i,r}$, for any $0\leq 
  i\leq [{k+m\over 2}]$. This follows from the above.
\end{pf}

If we introduce extra relations in $HFF^*_g$, we can expect that
elements of less degree may become zero acting on $\tilde{\cH}_r$. 
The next result deals with the case 
where we kill all elements $\q_i$, $1 \leq i \leq 2g$.

\begin{prop}
\label{prop:b1=0}
  We have that $\rb^{d+1}$ acts as zero on
  $\tilde{\cH}_r/(\seq{\q}{1}{2g})$,
  where $d=[{g-|r|-1 \over 2}]$.
\end{prop}

\begin{pf}
  From equation~\eqref{eqn:efe} one has that 
$$
  \rb^{d+1} =\sum_{i+j\geq d+1, j>0} c_{ij} \rb^i\g^{j},
$$
  on $\tilde{R}_{g,0,r}$, with $d=[{g-|r|-1 \over 2}]$. This means
  that $\rb^{d+1}$ sends $\tilde{\cH}_r$ into the ideal $(\seq{\q}{1}{2g})
  \subset \tilde{\cH}_r$, which is equivalent to the statement.
\end{pf}

Now we shall weaken the conditions in proposition~\ref{prop:b1=0}
in the form we need them to prove theorem~\ref{thm:C}. 
The following result is parallel to~\cite[proposition 6.12]{OS}. In
fact the proof of~\cite[proposition 6.12]{OS} can be done along 
the lines of the proof herein (see~\cite[section 8]{wang-yo}).

\begin{prop}
\label{prop:Szabo}
  Let $l=g-|r|-1$. Consider the ideal $\cI=(\seq{\q}{1}{l})$ in
$\tilde{\cH}_r$.
  Then any element $b \in\tilde{\AA}(\S)$ of degree bigger or equal
than $l+1$,
  sends $\tilde{\cH}_r$ into $\cI$.
\end{prop}

\begin{pf}
  It is enough to consider elements of degree $l+1$. A basic monomial
of 
  degree $l+1$ is of the form
  $$
     z=\rb^a \prod_{i\in I} (\q_i \q_{g+i}) \prod_{j \in J} \q_j
     \prod_{k \in K} \q_{g+k},
  $$
  where $I$, $J$, $K$ are disjoint subsets of $\{1,\ldots, g\}$ with
  $2a+2|I|+|J|+|K| = l+1$. If any $i\in I$ or $j\in J$ lies in 
  $\{1,\ldots, l\}$ then $z \in \cI$. Therefore we may suppose that
  $I$ and $J$ are disjoint with $\{1,\ldots, l\}$. 
  Moreover if $a=|I|=|J|=0$ then $z= \prod_{k \in K} \q_{g+k}\in 
  \L^{g-|r|}_0$ and by~\eqref{eqn:Hr2} then
  $z=0 \subset \cH_r$. So we also suppose
  $a+|I|+|J|>0$. Hence $|\{1,\ldots, l\} - K| \geq l -|K|
=2a+2|I|+|J|-1
  \geq |I|$.
  This means that we may associate (in an injective way) to every 
  $i\in I$ an $r(i)\in \{1,\ldots, l\} - K$,
  so that $\{ r(i) | i \in I\}$, $I$, $J$ and $K$
  are mutually disjoint. Therefore $z\in\cH_r$ is congruent (modulo
$\cI$) to
  $$
     z'=\rb^a \prod_{i\in I} (\q_i \q_{g+i}-\q_{r(i)} \q_{g+r(i)}) 
    \prod_{j \in J} \q_j \prod_{k \in K} \q_{g+k}.
  $$
  
  We conclude that $z' =\rb^a z_0$, with $z_0 \in \L_0^{l+1-2a}$.
  Now equation~\eqref{eqn:efe} tells us that
$$
\rb^a=\sum_{i+j\geq a,j>0} c_{ij}\rb^i\g^j,
$$
  acting on $\tilde{R}_{g,l+1-2a,r}$. So 
  we rewrite $z'$ as a polynomial with all homogeneous components of
  degree bigger or equal than $l+1$, and such that
  all its monomials contain powers 
  $\rb^i$, with $i<a$. By induction we are done, since when $a=0$,
  $z'=z_0 \in \L_0^{g-|r|}$, whose image is zero in $\cH_r$.
\end{pf}

\section{Proof of main results}
\label{sec:proofs}

\noindent {\em Proof of theorem~\ref{thm:0}\/}.
Suppose that $X$ is a $4$-manifold with $b^+>1$. In the first place
we shall reduce to the case of self-intersection $\S^2=0$ with
$\S$ representing an odd homology class. Suppose that $N=\S^2 >0$.
Take
the $N$-th blow-up of $X$ and call it $\tilde X=X \# N
\overline{\CP}^2$.
Let $\seq{E}{1}{N}$ denote the cohomology classes of the exceptional 
divisors, which are represented by embedded spheres of
self-intersection
$-1$. Consider the proper transform $\tilde \S$. This is an embedded
surface of genus $g$ in $\tilde X$ representing the homology class 
$\S -E_1 -\ldots -E_N$. It is obtained by tubing together $\S$ with
the exceptional spheres with reversed orientation inside $\tilde X$.
Then $\tilde \S^2=0$ and $\tilde\S$ represents an odd homology class,
since
if we take $w=E_1 \in H^2(X;\ZZ)$ then $\tilde \S \cdot w =1$.
On the other hand, the blow-up formula~\cite{bl} implies that the
order of 
finite type of $\tilde X$ is the same as that of $X$, so it is enough
to
prove the theorem for $\tilde X$ and $\tilde \S$.

Once we have done this reduction, choose $w \in H^2(X;\ZZ)$ with
$w \cdot \S \equiv 1\pmod 2$ (from~\cite{basic} we know that the
order 
of finite type does not depend on $w\in H^2(X;\ZZ)$). What we shall
prove
is the following: for any $b\in\tilde{\AA}(\S)$ with $d(b) \geq
2g-1$, we
have $D_X^w (b z )=0$, for any $z\in \AA(X)$. This is slightly
stronger
than the statement of the theorem.

We can suppose without loss of generality that $b$ is homogeneous,
i.e.
of the form $b=\wp^n\g_{i_1}\cdots \g_{i_m}$, where $d(b)=2n+m$.
Now consider a small tubular neighbourhood $A=\S \x D^2$ of $\S$
inside
$X$ and let $X_1$ be the complement of the interior of $A$. Therefore
$\bd X_1=Y=\Y$ and $X=X_1\cup_Y A$. Let $D \in H_2(X)$ with $D \cdot
\S=1$.
Represent $D$ by a $2$-cycle of the form $D=D_1 +\D$, where $D_1
\subset
X_1$ with $\bd D_1 =\SS^1$ and $\D=\point \x D^2 \subset A$. Take 
$z \in \AA(X_1)$. Then
$$
  \p^w(X_1, z e^{tD_1}) \in \whff
$$
may be decomposed according to components
$$
  \p^w(X_1, z e^{tD_1})_r \in \tilde{\cH}_r,
$$
with $-(g-1) \leq r \leq g-1$. With the notations of
section~\ref{sec:hff}, 
$\rb=\b +(-1)^{r+1}8$ in $\cH_r$. By corollary~\ref{cor:e},
$\rb^n \q_{i_1} \cdots \q_{i_m}$ is zero
acting on $\tilde{\cH}_r$, since it has degree
$d(b)=2n+m \geq 2g-1 \geq 2(g-|r|)-1$. Hence
$$
  0= {1\over 16^n} (\b^2-64)^n\q_{i_1} \cdots \q_{i_m}
  \p^w(X_1, z e^{tD_1})_r = \p^w(X_1,\wp^n\g_{i_1} \cdots \g_{i_m} 
   z e^{tD_1})_r,
$$
for any $r$. So $\p^w(X_1,b z e^{tD_1})=0$ and 
$$
  \Dws_X(b z e^{tD})=\la\p^w(X_1, bz e^{tD_1}),\p^w(A, e^{t\D})
\ra=0.
$$
In particular $D^w_X(b z e^{tD})=0$ for any $z \in \AA(X_1)$ and any 
$D \in H_2(X)$ with $D \cdot \S=1$. This implies the result
for a $4$-manifold with $b^+>1$. 
\hfill $\Box$

\vspace{5mm}
\noindent {\em Proof of theorem~\ref{thm:A}\/}.
First we shall reduce to the case of $\S^2=0$ and $\S$ representing
an odd homology class.
Let $X$ be a $4$-manifold with $b^+>1$ and consider its blow-up
$\tilde X=X\# \overline{\CP}^2$ at one point, where $E$ stands for 
the cohomology class of the exceptional divisor.
Now if $b \in \tilde{\AA}(X)=\tilde{\AA}(\tilde X)$ and $K$ is
a basic class for $D^w_X(b\,\bullet)$, then there exists $z \in
\AA(X)$ such
that $D^w_X(bze^{tD+\l x})=e^{Q(tD)/2 +2\l + K\cdot tD}$, for all 
$D \in H_2(X)$. The blow-up formula~\cite{bl} tells us that  
$D^w_{\tilde X}(bze^{tD+\l x})=e^{Q_{\tilde X}(tD)/2 +2\l + K \cdot
tD} 
\cosh (E\cdot tD)$, for $D\in H_2(\tilde X)$ 
(compare~\cite[proposition 16]{basic}). Therefore $K \pm E$
are basic classes for $D^w_{\tilde X}(b\,\bullet)$. If $\S$ is an
embedded
surface in $X$ of genus $g$ with self-intersection $\S^2 >0$ then its
proper transform
$\tilde \S$ is an embedded surface in $\tilde X$ of genus $g$,
$\tilde \S^2=\S^2-1$ and representing the homology class $\S-E$. Put
$\e$ for the sign of
$K\cdot\S$. Then $K-\e E$ is a basic class for $\tilde X$ and 
$$
  |(K - \e E )\cdot \tilde{\S} | +\tilde \S^2 = |K \cdot \S| + \S^2.
$$
So if we prove the theorem for $\tilde X$ then it will be proved for
$X$.
This means that we can reduce to the case of self-intersection
$\S^2=0$
and $\S$ representing an odd element in homology.

As in the proof of theorem~\ref{thm:0}, we choose $w \in H^2(X;\ZZ)$
with
$w\cdot \S \equiv 1\pmod 2$ (since $\S$ represents an odd homology
class).
Also we can suppose $b=\wp^n\g_{i_1}\cdots \g_{i_m}$, where
$d(b)=2n+m$.
Consider a small tubular neighbourhood $A=\S \x D^2$ of $\S$ inside
$X$ and let $X_1$ be the complement of the interior of $A$. Therefore
$\bd X_1=Y=\Y$ and $X=X_1\cup_Y A$. Let $D \in H_2(X)$ with $D \cdot
\S=1$.
Represent $D$ by a $2$-cycle of the form $D=D_1 +\D$, where $D_1
\subset
X_1$ with $\bd D_1 =\SS^1$ and $\D=\point \x D^2 \subset A$.
Fix an homogeneous $z_0 \in \L^* H_1(X) \subset \AA(X)$. 
Then $\p^w(X_1, b z_0 e^{tD_1}) \in \whff$ and
\begin{equation}
  \Dws_X(b z_0 e^{tD +\l x +s\S})=\la\p^w(X_1, b z_0
e^{tD_1}),\p^w(A, e^{t\D
  +\l x+s\S}) \ra.
\label{eqn:k}
\end{equation}
Decomposing according to components $\whff=
\oplus \tilde{\cH}_r$, we have
$$
  \Dws_X(b z_0 e^{tD +\l x +s\S})=\sum_{r=-(g-1)}^{g-1}
  \la\p^w(X_1, b z_0 e^{tD_1})_r,\p^w(A, e^{t\D +\l x+s\S})_r \ra.
$$
Put $d_0=-w^2- {3\over 2}(1-b_1+ b^+)$, so
$$
  D^w_X(b z_0 e^{tD +\l x +s\S})= {1 \over 2}\left( \Dws_X(b z_0
e^{tD +\l x +s\S})
  +\ima^{-d_0+d(b)+d(z_0)}\Dws_X(b z_0 e^{\ima tD -\l x +\ima s\S})
\right)=
$$
$$
 ={1\over 2} \Big( \sum_{r=-(g-1)}^{g-1}
  \la\p^w(X_1, b z_0 e^{tD_1})_r,\p^w(A, e^{t\D +\l x+s\S})_r \ra + 
$$
\begin{equation}
\label{eqn:33}
  + \ima^{-d_0+d(b)+d(z_0)}\sum_{r=-(g-1)}^{g-1}
  \la\p^w(X_1, b z_0 e^{\ima tD_1})_r,\p^w(A, e^{\ima t\D -\l x+\ima
s\S})_r \ra
  \Big).
\end{equation}
By the definition of $\cH_r$ in~\eqref{eqn:Hr} through the
eigenvalues of $\a$, we have that $\p^w(A, e^{t\D +\l x+s\S})_r$ 
is killed by some product of differential operators of the
form ${\bd \over \bd s}-(2r+f_i(t))$, $f_i(t)\in t\Ct$, 
if $r$ is odd, and by some product of differential operators 
of the form ${\bd \over \bd s}-(2r\ima+g_i(t))$, 
$g_i(t)\in t\Ct$, if $r$ is even. Also $\p^w(A, 
e^{\ima t\D -\l x+\ima s\S})_r$ is killed by some product of 
differential operators ${\bd \over \bd s} -(2r\ima+f_i(\ima t)\ima)$ 
if $r$ is odd, and by some product of differential operators 
${\bd \over \bd s} -(-2r+g_i(\ima t)\ima)$ if $r$ is even. 
So the existence of an exponential term of the form
$e^{(2r+f(t))s}$, $f(t)\in t\Ct$, 
in~\eqref{eqn:33} is equivalent to the non-vanishing of
the summand corresponding to $r$.

Put $2r=K\cdot \S$ (it is an even number since $K$ is an integral
lift 
of $w_2(X)$ by~\cite[theorem 6]{basic}). 
That $K$ is a basic class for $D^w_X(b\,\bullet)$ means that there is 
$z \in \AA(X)$ (which we may suppose to be of the form $z=z_0 z_1$,
with
$z_0\in \L^* H_1(X)$ homogeneous and $z_1 \in \Sym^* (H_0(X)\oplus
H_2(X))$)
with $D^w_X(b z_0z_1 e^{tD'+\l x})=e^{Q(tD')/2+2\l+K\cdot tD'}$,
for any $D'\in H_2(X)$. So
$$
 D^w_X(bz_0z_1e^{tD+\l x+s\S})=e^{Q(tD)/2+2\l+K\cdot tD+(2r+t)s}.
$$
Hence an exponential term of the form $e^{(2r+t)s}$ appears
in~\eqref{eqn:33}
and thus $\p^w(X_1, bz_0 e^{tD_1})_r \neq 0$. As $b=\wp^n \g_{i_1}
\cdots \g_{i_m}$,
we deduce that $\rb^n \g_{i_1}\cdots \g_{i_m}$ acts 
as non-zero on $\tilde{\cH}_r$ (note that
$\b+(-1)^r8$ is an isomorphism in $\tilde{\cH}_r$).
By corollary~\ref{cor:e}, $d(b)=2n+m \leq 2(g-|r|-1)$. So $|2r|+d(b)
\leq 2(g-1)$.
\hfill $\Box$

\vspace{5mm}
\noindent {\em Proof of theorem~\ref{thm:B}\/}.
In principle theorem~\ref{thm:B} follows from the more general
theorem~\ref{thm:C},
but we shall give an independent proof because of its simplicity.
It is similar to the proof of theorem~\ref{thm:A}. We reduce to
consider
the case $\S^2=0$ and $\S$ representing an odd element in homology.
Now by definition~\ref{def:d(K)}, and since $X$ has $b_1=0$, we have
that
$K$ is a basic class for $D^w_X(b \, \bullet)$ where $b=\wp^n$ and
$2n=d(K)$.

We perform the splitting $X=X_1 \cup_Y A$ as in 
the proof of theorem~\ref{thm:A}. Then 
$\p^w(X_1, e^{tD_1}) \in \whff$ lies in the kernels of all
$\q_i$, $1\leq i \leq 2g$, since $X$ has $b_1=0$. 
By proposition~\ref{prop:b1=0}, $\rb^{d+1}$ sends $\tilde{\cH}_r$
to the ideal $(\seq{\q}{1}{2g}) \subset \tilde{\cH}_r$,  
where $d=[{g-|r|-1 \over 2}]$. In particular, 
$\rb^{d+1}\p^w(X_1, e^{tD_1})_r=0$.

Put $2r =K \cdot \S$. Again 
$K$ being a basic class for $D^w_X(b \, \bullet)$ implies
that $\p^w(X_1, b e^{tD_1})_r \neq 0$. Therefore $\rb^n \neq 0$ and 
hence $n \leq d =[{g-|r|-1 \over 2}]$. So $2n \leq g-|r|-1$,
i.e.\ $|2r| +2d(K) \leq 2(g-1)$. \hfill $\Box$

\vspace{5mm}
\noindent {\em Proof of theorem~\ref{thm:C}\/}.
As above we reduce
to the case $\S^2=0$ and $\S$ representing an odd element in
homology.
We may suppose without loss of generality that $b=\wp^n \g_{i_1}
\cdots
\g_{i_m}$ with $d(b)=2n+m$. Keeping the notations of the proof of 
theorem~\ref{thm:B}, we have that
$$
  \p^w(X_1,b e^{tD_1})_r \in \tilde{\cH}_r
$$
is non-zero, where $2r =K\cdot D$. 

If $l+1 \leq g-1-|r|$ then obviously $|2r| +2d(b) \leq 2(g-1)$ is
true.
Otherwise $g-1-|r| \leq l$. Then $\p^w(X_1,b e^{tD_1})_r$
lives in the kernels of $\seq{\q}{1}{l}$, since $i_*(\g_j)=0 \in
H_1(X)$
for $j=1,\ldots, l$. Therefore we conclude that 
$\rb^n \q_{i_1} \cdots \q_{i_m} \neq 0$ acting on
$\cH_r/(\seq{\q}{1}{g-|r|-1})$. By proposition~\ref{prop:Szabo},
we have that $d(b) \leq g-|r|-1$, i.e.\  $|2r| +2d(b) \leq 2(g-1)$.
\hfill $\Box$

\noindent {\em Acknowledgements:\/} I am grateful to the Department
of
Mathematics of Universidad Au\-t\'o\-noma de Madrid for their hospitality
during my stay in the first semester of the year 98/99. Also thanks
to
Peter Ozsv\'ath and Zolt\'an Szab\'o for useful correspondence.


\begin{thebibliography}{10}
\bibitem{basic}
V. Mu\~noz,
Basic classes for four-manifolds not of simple type,
{\tt math.DG/9811089}.
\end{thebibliography}
\end{document}